\documentclass[a4paper]{article}

\usepackage{latexsym}
\usepackage{amssymb}
\usepackage{amsthm}
\usepackage{amsmath}
\usepackage{amssymb}

\usepackage{
  amstext,
  amsbsy,
  amsopn,
  amsfonts
  }

\newcommand{\p}{\partial}
\newcommand{\eps}{\varepsilon}
\newcommand{\vphi}{\varphi}

\newcommand{\dd}{\hspace{1pt}{\rm d}\hspace{0.0pt}}

\newcommand{\R}{\mathbb R}
\newcommand{\N}{\mathbb N}
\newcommand{\C}{\mathbb C}

\newcommand{\cC}{{\mathcal C}}
\newcommand{\cD}{{\mathcal D}}
\newcommand{\cE}{{\mathcal E}}

\newcommand{\cG}{{\mathcal G}}
\newcommand{\cL}{{\mathcal L}}
\newcommand{\cN}{{\mathcal N}}
\newcommand{\cO}{{\mathcal O}}
\newcommand{\cP}{{\mathcal P}}

\newcommand{\cS}{{\mathcal S}}

\newcommand{\rL}{\ensuremath{\mathrm{L}}}
\newcommand{\Cinf}{\ensuremath{{\mathcal C}^\infty}}

\newcommand{\EM}{\ensuremath{{\cE}_{\mathrm{M}}}}
\newcommand{\Ginf}{\ensuremath{\cG^\infty}}

\newcommand{\supp}{\mathop{\mathrm{supp}}}

\newcommand{\hp}{\ensuremath{\R\times[0,\infty)}}

\newcommand{\comp}{^{\sf C}}
\newcommand{\eins}{\mathbf 1}

\newtheorem{thm}{Theorem}

\newtheorem{prop}[thm]{Proposition}
\newtheorem{defn}[thm]{Definition}

\newtheorem{rem}[thm]{Remark}
\newtheorem{ex}[thm]{Example}

\begin{document}

\title{Colombeau solutions to nonlinear wave equations}

\vspace{5mm}
\author{Michael Oberguggenberger\thanks{Partially supported by FWF(Austria), grant Y237 and by FAPESP (Brasil)}\\
{\normalsize Institut f\"ur Grundlagen der Baungenieurwissenschaften}\\
{\normalsize Universit\"at Innsbruck, A-6020 Innsbruck, Austria}\\
{\normalsize E-mail: michael.oberguggenberger@uibk.ac.at}}
\date{}
\maketitle
\begin{abstract}
This paper is a tutorial that demonstrates various methods from
the Colombeau theory of generalized functions in the context of
semilinear wave equations. The Colombeau generalized functions
constitute differential algebras that contain the space of
distributions. We solve the 1D-semilinear wave equation in these
algebras, show how delta waves can be computed and then turn to
linear and nonlinear regularity theory.
\vspace{8pt}
\end{abstract}

%
%
\section{Introduction}
\label{Sec:intro}
%
%
This tutorial serves the purpose of demonstrating how the theory of Colombeau algebras
can be used to solve and study nonlinear partial differential equations that do not
have solutions in the sense of distributions.
It starts out with existence results obtained in the eighties of the twentieth century
and leads up to recent regularity results.

For the purpose of exposition we focus on a model problem,
the one-dimensional semilinear wave equation with singular initial data
\begin{equation}\label{eq:semilinearwave}
\begin{array}{l}
\p_t^2 u(x,t) - \p_x^2 u(x,t) = f(u(x,t)) + h(x,t),\quad (x,t)\in\R^2,\\
u(x,0) = a(x),\ \p_tu(x,0) = b(x),\quad x\in\R,
\end{array}
\end{equation}
where $f$ is a smooth, polynomially bounded function and $a, b$ and $h$ are Colom\-beau
generalized functions on the real line and on the $(x,t)$-plane, respectively.
The main questions to be addressed are the following:

(a) existence and uniqueness of generalized solutions in Colombeau algebras;\\
(b) limiting behavior of the representatives when the data are distributions;\\
(c) regularity of generalized solutions.

Concerning question (a), a short survey of existence and uniqueness results
can be found in \cite{MOWCNA}. The existence of Colombeau solutions to semilinear hyperbolic systems
was one of the early applications of the theory \cite{Bia88, MO87}. For linear hyperbolic systems
with non-smooth (Colombeau generalized) coefficients, existence and
uniqueness was established in \cite{MO89}, for symmetric hyperbolic systems in higher space dimensions in
\cite{LafMO90} and for hyperbolic pseudodifferential systems with Colombeau symbols in
\cite{Hoer03}. Notably for semilinear hyperbolic systems, question (b)
has been answered in many cases, involving the notion of delta-waves
\cite{ColMO90, Gram91, MO86, MO91, MOWang94a, MOWang94b, RauReed87}.
Finally, regularity theory for Colombeau solutions is now based on the subalgebra
$\Ginf$ of regular Colombeau functions and currently an active area of research,
making use of pseudodifferential and microlocal techniques
\cite{Gar04, GarGramMO03, HoerHoop01, NedPilSca98}. In particular, the
propagation of the $\Ginf$-wave front set in linear systems with
Colombeau coefficients is a theme of recent investigations \cite{GarHoer04}.
To date, only few results are available about regularity of Colombeau solutions in the nonlinear case.

The plan of the tutorial is as follows.
After recalling the required notions from Colombeau theory in Section \ref{Sec:notation},
we shall prove existence and uniqueness of a solution $u$ belonging to the
Colombeau algebra on $\R^2$ in Section \ref{Sec:exis} and state various preliminary
regularity properties.
Section \ref{Sec:sing} is devoted to computing the associated distribution
(the distributional limit of the representing nets) when the initial data are
delta functions or derivatives of delta functions. In Section \ref{Sec:linregu}
we turn to regularity theory. We recall the $\Ginf$-regularity result
for the linear case and show that it fails in the nonlinear case.
The tutorial finishes with a recent result on propagation of regularity in the nonlinear
case in Section \ref{Sec:linregu}.

For simplicity of presentation, we restrict our attention to the one-dimensional
case and Lipschitz-continuous nonlinearity $f$. At the appropriate places of the
tutorial, we will indicate what is known about the non-Lipschitz and the
higher dimensional case. With Lipschitz-continuous $f$, in particular, the
solution to (\ref{eq:semilinearwave}) exists globally in space and time. Therefore, we
need not enter the discussion of domains of existence nor use other versions
of Colombeau algebras to accommodate energy estimates when $f$ grows polynomially.
%
%
\section{Notation}
\label{Sec:notation}
%
%
The paper is placed in the framework of algebras of generalized functions
introduced by Colombeau in \cite{Col84, Col85}.
We shall fix the notation and introduce a number of known as well as new classes
of generalized functions here. For more details, see \cite{GKOS01}.

Let $\Omega$ be an open subset of $\R^n$. The basic objects of the theory as we
use it are families
$(u_\eps)_{\eps \in (0,1]}$ of smooth functions $u_\eps \in \Cinf(\Omega)$ for
$0 < \eps \leq 1$.
We single out the following subalgebras:

{\em Moderate families}, denoted by $\EM(\Omega)$, are defined by the property:
\begin{equation}
  \forall K \Subset \Omega\,\forall \alpha \in \N_0^n\,
    \exists p \geq 0:\;\sup_{x\in K} |\p^\alpha u_\eps(x)|
       = \cO(\eps^{-p})\ \rm{as}\ \eps \to 0.
     \label{eq:mofu}
\end{equation}
{\em Null families}, denoted by $\cN(\Omega)$, are defined by the property:
\begin{equation}
  \forall K \Subset \Omega\,\forall \alpha \in \N_0^n\,
     \forall q \geq 0:\;\sup_{x\in K} |\p^\alpha u_\eps(x)|
        = \cO(\eps^q)\ \rm{as}\ \eps \to 0.
      \label{eq:nufu}
\end{equation}
Thus moderate families satisfy a locally uniform polynomial estimate as $\eps \to 0$,
together with all derivatives, while null functionals
vanish faster than any power of $\eps$ in the same situation. The null families
form a differential ideal in the collection of moderate families. The {\em Colombeau
algebra} is the factor algebra
\[
   \cG(\Omega) = \EM(\Omega)/\cN(\Omega).
\]
The algebra $\cG(\Omega)$ just defined coincides with the {\em special Colombeau
algebra} in \cite[Def. 1.2.2]{GKOS01}, where the notation $\cG^s(\Omega)$
has been employed. It was called the {\em simplified Colombeau algebra} in \cite{Biabook}.

Families $(r_\eps)_{\eps\in(0,1]}$ of complex numbers such that $|r_\eps| = \cO(\eps^{-p})$
as $\eps \to 0$ for some $p \geq 0$ are called {\em moderate}, those for which
$|r_\eps| = \cO(\eps^q)$ for every $q \geq 0$ are termed {\em negligible}.
The ring $\widetilde{\C}$ of Colombeau generalized numbers is obtained by factoring moderate families
of complex numbers with respect to negligible families. When $\Omega$ is connected, $\widetilde{\C}$
coincides with the {\em ring of constants} in the differential algebra $\cG(\Omega)$.

The restriction of an element $u \in \cG(\R^2)$
to the line $\{t=0\}$
is defined on representatives by
\[
   u\vert_{\{t=0\}}\ = \ \mbox{class of}\ (u_\eps(\cdot,0))_{\eps \in (0,1]}.
\]
Similarly, restrictions of the elements of $\cG(\Omega)$ to open subsets of $\Omega$ are defined
on representatives. One can see that $\Omega \to \cG(\Omega)$ is a sheaf of differential
algebras on $\R^n$.
The space of compactly supported distributions is imbedded in $\cG(\Omega)$ by convolution:
\begin{equation} \label{eq:imbedding}
   \iota:\cE'(\Omega) \to \cG(\Omega),\
     \iota(w)\ = \ \mbox{class of}\ (w \ast (\varphi_\eps)\vert_\Omega)_{\eps \in (0,1]},
\end{equation}
where
\begin{equation} \label{eq:molli}
   \varphi_\eps(x) = \eps^{-n}\varphi\left(x/\eps\right)
\end{equation}
is obtained by scaling a fixed test function $\varphi \in \cS(\R^n)$ of integral one
with all moments vanishing. By the sheaf property, this can be extended in a unique
way to an imbedding of the space of distributions $\cD'(\Omega)$. In the case
$\Omega = \R^n$, the imbedding is given by the explicit formula
\begin{equation}
   \iota:\cD'(\Omega) \to \cG(\R^n),\
     \iota(w)\ = \ \mbox{class of}\ (w \ast (\chi\varphi_\eps))_{\eps \in (0,1]},
\end{equation}
where $\chi$ is some compactly supported smooth function identically equal to one in
a neighborhood of zero.  We refer to \cite[Sec. 1.2]{GKOS01} for further
explicit expressions in general and in special cases.

One of the main features of the Colombeau construction is the fact that this imbedding
renders $\Cinf(\Omega)$ a faithful subalgebra. In fact, given $f \in \Cinf(\Omega)$,
one can define a corresponding element of $\cG(\Omega)$ by the constant imbedding
$\sigma(f)\ = \ \mbox{class of}\ [(\eps,x) \to f(x)]$.
Then the important equality $\iota(f) = \sigma(f)$ holds in $\cG(\Omega)$.

If $u \in\cG(\Omega)$ and $f$ is a smooth function
which is of at most polynomial growth at
infinity, together with all its derivatives, the superposition $f(u)$ is a well-defined
element of $\cG(\Omega)$.

The algebra $\cG(\Omega)$ can be equipped with a topology that turns it into a
topological ring and a topological module over $\widetilde{\C}$, and into a complete
ultrametric space. Let $K\Subset\Omega$, $m\in\N_0$ and $p \geq 0$.
Then $V(K,m,p)\subset\cG(\Omega)$ is defined as the collection of elements of $u\in\cG(\Omega)$
with a representative $u_\eps$ such that
\[
   \sup_{x\in K} \sup_{|\alpha| \leq m}|\p^\alpha u_\eps(x)|
       = \cO(\eps^{-p})\ \rm{as}\ \eps \to 0.
\]
The sets $V(K,m,p)\subset\cG(\Omega)$ define a base of neighborhoods of zero for a topology
with the properties mentioned above. This topology was introduced
in \cite{Biabook}; its vast potential was discovered by \cite{Scarpalezos:00} who coined the
term {\em sharp topology}. Important further developments are due to \cite{AragonaJuriaans,Garetto:05a}.

We need a couple of further notions from the theory of Colombeau generalized functions.
Regularity theory  is based on
the subalgebra $\cG^\infty(\Omega)$ of {\em regular
generalized functions} in $\cG(\Omega)$. It is defined by those elements which have
a representative satisfying
\begin{eqnarray}
  \forall K \Subset \Omega\,\exists p \geq 0\,\forall \alpha \in \N_0^n:\;
  \sup_{x\in K} |\p^\alpha u_\eps(x)| = \cO(\eps^{-p})\
     \ \rm{as}\ \eps \to 0. \label{eq:regufu}
\end{eqnarray}
Observe the change of quantifiers with respect to formula (\ref{eq:mofu}); locally, all derivatives
of a regular generalized function have the same order of growth in $\eps > 0$. One has
that (see \cite[Thm. 5.2]{MObook})
\[
  \cG^\infty(\Omega) \cap \cD'(\Omega) = \Cinf(\Omega).
\]
For the purpose of describing the regularity of Colombeau generalized functions,
$\cG^\infty(\Omega)$ plays the same role as $\Cinf(\Omega)$ does in the setting of distributions.
However, $\Ginf(\Omega)$ is not invariant under nonlinear maps. For this reason, various
subspaces of $\cG(\Omega)$ have been introduced to measuring regularity in the nonlinear case
(see e.g. \cite{MONoviSad}). We shall make use of just one of them, based on the notion of
subsheaf regularity introduced in \cite{Marti:99}.
\begin{defn}\label{def:subsheaf}
$\cL_\cG^1(\Omega)$ is the space of the elements of $\cG(\Omega)$ with a representative
$(u_\eps)_{\eps \in (0,1]}$ such that $\lim_{\eps \to 0} u_\eps$ exists in $\rL_{\mathrm{loc}}^1(\Omega)$.
\end{defn}
\begin{prop}\label{prop:superpossub}
$f\big(\cL_\cG^1(\Omega)\big) \subset \cL_\cG^1(\Omega)$ for every smooth function $f$ all whose derivatives
grow at most polynomially at infinity and which is Lipschitz continuous with a global Lipschitz constant.
\end{prop}
{\it Proof:} The polynomial bounds guarantee that $f(u)$ is a well-defined element of $\cG(\Omega)$.
Assume that $u_\eps$ converges to an element $w \in \rL^1(K)$ on some compact set $K$. The estimate
$|f(u_\eps) - f(w)| \leq \mathrm{Lip}_f|u_\eps - w|$ shows that $f(u_\eps)$ converges to $f(w)$, as
desired.\qed

The kind of regularity that is encapsulated in the subspace above is described by the following
rather obvious assertion:
\[
   \cL_\cG^1(\Omega) \cap \cD'(\Omega) = \rL_{\mathrm{loc}}^1(\Omega).
\]
We end this section by recalling the {\em association relation} on the Colombeau
algebra $\cG(\Omega)$. It identifies elements of $\cG(\Omega)$ if they
coincide in the weak limit. That is, $u, v \in \cG(\Omega)$ are called associated,
$u \approx v$, if
$ \lim_{\eps \to 0} \int\big(u_\eps(x) - v_\eps(x)\big) \psi(x)\,dx = 0 $
for all test functions $\psi \in \cD(\Omega)$. We shall also say that
$u$ is associated with a distribution $w$ if $u_\eps \to w$
in the sense of distributions as $\eps\to 0$.
%
%
\section{Existence/uniqueness of generalized solutions}
\label{Sec:exis}
%
%
This section is devoted to solving the semilinear wave equation (\ref{eq:semilinearwave})
in the Colombeau algebra $\cG(\R^2)$. Recall first that if $w$ is a classical
solution of the linear wave equation
\begin{equation}\label{eq:linearwave}
\begin{array}{l}
\p_t^2 w(x,t) - \p_x^2 w(x,t) = h(x,t),\quad (x,t)\in\R^2,\\
w(x,0) = a(x),\ \p_tw(x,0) = b(x),\quad x\in\R,
\end{array}
\end{equation}
then it solves the integral equation
\begin{equation}\label{eq:linearwaveintegral}
  w(x,t) = \frac{1}{2}\big(a(x-t) + a(x+t)\big) + \frac{1}{2}\int_{x-t}^{x+t}b(y) \dd y
    + \frac{1}{2} \int_0^t\int_{x-t+s}^{x+t-s} h(y,s)\dd y\dd s.
\end{equation}
Let $K_0 = [-\kappa,\kappa]$ be a compact interval. For $0\leq T \leq \kappa$,
the trapezoidal region $K_T$ is defined by
\begin{equation} \label{eq:ItKT}
  K_T = \{(x,t)\in \R^2 : 0\leq t\leq T, |x| \leq \kappa - t\}.
\end{equation}
Using (\ref{eq:linearwaveintegral}), the following estimate is
easily deduced ($0\leq t \leq T \leq \kappa$):
\begin{eqnarray}
\|w\|_{\rL^\infty(K_T)} \leq \|a\|_{\rL^\infty(K_0)} + T\|b\|_{\rL^\infty(K_0)}
    + T \int_0^T \|h\|_{\rL^\infty(K_s)}\dd s.     \label{eq:Linftyestimate}
\end{eqnarray}
For later reference, we note that estimate (\ref{eq:Linftyestimate}) holds with the
$\rL^1$-norm in place of the $\rL^\infty$-norm as well.
We now turn to the semilinear wave equation (\ref{eq:semilinearwave}). We assume
throughout that $u \to f(u)$ is a smooth function all whose derivatives are of at most
polynomial growth as $|u|\to\infty$, that $f$ satisfies a global Lipschitz
estimate (i.e., has a bounded first derivative) and that $f(0) = 0$.
\begin{prop} \label{prop:exis}
Assume that the function $f$ is as described above. Let $a,b \in \cG(\R)$ and $h \in \cG(\R^2)$.
Then problem (\ref{eq:semilinearwave}) has a unique solution $u \in \cG(\R^2)$.
The solution depends continuously on the data with respect to the sharp topology.
\end{prop}
{\it Proof:} To prove the existence of a solution, take representatives
$a_\eps, b_\eps, h_\eps$ of $a,b,h$, respectively, and let
$u_\eps \in \Cinf(\R^2)$ be the unique classical solution to the semilinear wave
equation with regularized data:
\begin{equation}\label{eq:waverepresentatives}
\begin{array}{lcl}
\p_t^2 u_\eps - \p_x^2 u_\eps = f(u_\eps) + h_\eps &\quad \mbox{on}& \R^2,\\
u_\eps(\cdot,0) = a_\eps,\ \p_tu_\eps(\cdot,0) = b_\eps &\quad \mbox{on}& \R.
\end{array}
\end{equation}
The classical solution $u_\eps$ to (\ref{eq:waverepresentatives}) is constructed
by rewriting (\ref{eq:waverepresentatives}) as an integral equation and invoking a
fixed point argument (this involves applying estimate (\ref{eq:Linftyestimate}) successively to
all derivatives). If we show that the net $(u_\eps)_{\eps\in(0,1]}$ belongs to
$\EM(\R^2)$, its equivalence class in $\cG(\R^2)$ will be a solution.
We shall do the proof only for the upper half-plane $(x,t)\in\hp$; the arguments
for $t\leq 0$ are similar.
To show that the zero-th derivative of $u_\eps$ satisfies the estimate (\ref{eq:mofu}),
we take a region $K_T$ and invoke
inequality (\ref{eq:Linftyestimate}) to see that
\begin{equation}\label{eq:decisiveestimate}
     \|u_\eps\|_{\rL^\infty(K_t)} \leq \|a_\eps\|_{\rL^\infty(K_0)} +
      T\|b_\eps\|_{\rL^\infty(K_0)} \vspace{4pt}\\
     \displaystyle +\ T \int_0^t \|f(u_\eps) + h_\eps\|_{\rL^\infty(K_s)}\dd s.
  \end{equation}
The last term on the right hand side of (\ref{eq:decisiveestimate}) is estimated by
\begin{equation}\label{eq:decisiveestimatelast}
     T \int_0^t \|f'\|_{\rL^\infty(\R)}
       \|u_\eps\|_{\rL^\infty(K_s)}\dd s
           + T^2\|h_\eps\|_{\rL^\infty(K_T)}
\end{equation}
Using that each of the terms involving
$a_\eps, b_\eps, h_\eps$ is of order $\cO(\eps^{-p})$ for some $p$, we infer
from Gronwall's inequality that the same is true of $\|u_\eps\|_{\rL^\infty(K_t)}$
for $0 \leq t \leq T$. Thus $u_\eps$ is moderate on the region $K_T$, that is, it
satisfies the estimate (\ref{eq:mofu}) there. To get the estimates for the higher order derivatives,
one just differentiates the equation and employs the same arguments inductively, using
that the lower order terms are already known to be moderate from the previous steps.

To prove uniqueness, we consider representatives $u_\eps, v_\eps \in \EM[\R^2)$ of two solutions
$u$ and $v$. Their difference satisfies
\begin{equation*}
\begin{array}{l}
\p_t^2 (u_\eps - v_\eps) - \p_x^2 (u_\eps - v_\eps)
   = \left(f(u_\eps)-f(v_\eps)\right) + n_\eps,\\
     (u_\eps(\cdot,0)- v_\eps(\cdot,0))= n_{0\eps},\
      \p_t(u_\eps(\cdot,0) - v_\eps(\cdot,0)) = n_{1\eps}
\end{array}
\end{equation*}
for certain null elements $n_\eps, n_{0\eps}, n_{1\eps}$.
Thus $u_\eps - v_\eps$ satisfies an estimate of the form (\ref{eq:decisiveestimate}),
but with the null elements $n_\eps, n_{0\eps}, n_{1\eps}$ replacing
$a_\eps, b_\eps,h_\eps$ there. This implies as above that the $\rL^\infty$-norm
of $u_\eps - v_\eps$ on $K_T$ is of order $\cO(\eps^q)$ for every $q\geq 0$.
By \cite[Thm. 1.2.3]{GKOS01}, the null estimate (\ref{eq:nufu}) on $u_\eps - v_\eps$ suffices
to have null estimates on all derivatives. Thus $u = v$ in $\cG(\R^2)$.

The proof of continuous dependence on the data $a,b,h$ follows exactly the same lines.
\qed

Having established a general existence- and uniqueness result, the question arises
what else can be said about the solution, when more is known about the data. A number
of results on the qualitative properties are indeed available. We begin with the most
basic result on $\cC^\infty$-smoothness of the solution when the data are smooth.
\begin{prop}
Assume that $a,b$ belong to $\Cinf(\R)$ and $h$ is in $\Cinf(\R^2)$. Then the generalized solution
$u\in \cG(\R^2)$ to problem (\ref{eq:semilinearwave}) coincides with the classical solution
$w \in \Cinf(\R^2)$, that is, $u = \iota(w)$ in $\cG(\R^2)$.
\end{prop}
{\em Proof:} The imbeddings $\iota$ and $\sigma$ coincide on
$\cC^\infty(\R)$ and on $\cC^\infty(\R^2)$. Thus we may represent the data
by $a_\eps \equiv a$, $b_\eps \equiv b$ and $h_\eps \equiv h$. The proof of
Prop.\;\ref{prop:exis} shows that $u_\eps \equiv w$ is a representative of the
generalized solution, and this means that $u = \iota(w)$ in $\cG(\R^2)$. \qed
\begin{prop}
Assume that the data $a,b$ and $h$ are continuous functions and let $w \in \cC(\R^2)$
be the corresponding continuous (weak) solution. Then the generalized solution $u\in \cG(\R^2)$
is associated with $w$, that is, $u \approx \iota(w)$.
\end{prop}
{\em Proof:} Let $a_\eps = a \ast \varphi_\eps$ and similarly for $b_\eps$ and $h_\eps$ with appropriate
one- or two-dimensional mollifiers $\varphi_\eps$. On the one hand, the families $a_\eps, b_\eps, h_\eps$
define representatives of the generalized functions $\iota(a),\iota(b),\iota(h)$. On the other hand,
they converge locally uniformly to the original continuous functions $a,b,h$.
The estimates (\ref{eq:decisiveestimate}) and (\ref{eq:decisiveestimatelast}), applied to
$u_\eps - w$, $a_\eps - a$ etc. together with
Gronwall's inequality, show that the solution $u_\eps$ converge to $w$ with respect the sup-norm on compact sets.
In particular, $u$ is associated with $w$. \qed

Note that in general $u \neq \iota(w)$ as elements of $\cG(\R^{n+1})$. As a rule, equality
of continuous and generalized solutions does not hold in any algebra of generalized functions.
This is similar to Schwartz' impossibility result and has been worked out in \cite{MO89rough}.
When the data are distributions, there may be no meaning for a distributional solution,
in general. Yet the solution in $\cG(\R^2)$ may still be associated with a distribution.
We will turn to an incident of such a situation in the next section.

Proposition \ref{prop:exis} is a model result. In fact, for vanishing
driving term $h \equiv 0$, it is a special case of \cite{Bia88} and of \cite{MO87};
for Lipschitz continuous, smooth $f$, existence and uniqueness
of a solution in $\cG(\R^2)$ can be proven in space dimensions
$n = 1,2,3$, see e.g.\ \cite{MOWCNA}. If $f$ is not Lipschitz, but of polynomial growth, energy estimates
can be used to construct solutions in the Colombeau algebra $\cG_{2,2}(\R^n\times[0,\infty))$
introduced in \cite{BiaMO92b}. As in the classical case, the growth type
of $f$ is connected with the space dimension; the cases $1 \leq n \leq 6$ have been
treated in \cite{NedMOPil03}; for $n=3$ see also \cite{Col85}.
%
%
\section{Delta waves}
\label{Sec:sing}
%
%
In this section we investigate the behavior of the solution $u \in\cG(\R^2)$ to
the semilinear wave equation (\ref{eq:semilinearwave}) when the initial data
$a,b$ are given by distributions. In certain cases, the smooth solutions
$u_\eps$ with regularized initial data $a_\eps, b_\eps$ converge to a distributional
limit. When it exists, this limit is called a {\em delta wave}. Stated equivalently,
the generalized solution $u\in \cG(\R^2)$ with data $\iota(a),\iota(b)$ admits an
associated distribution. To be sure, this distribution is not a solution of the
original nonlinear wave equation in any sense. It rather has to be viewed as describing
the qualitative properties of the Colombeau solution.

For simplicity of exposition, we will assume here that the right hand side $h$ is zero
(see e.g. \cite{MOWCNA} for a result with nonzero $h$) and will exhibit the model
case of bounded nonlinearity $f$ and initial data with discrete support.
Thus $a, b \in \cD'(\R)$ will be assumed to be Dirac measures or their derivatives situated at a finite set
of points. Consider the weak solution $v \in \cD'(\R^{2})$ to the linear problem
\begin{equation}\label{eq:classicaldeltalinear}
\begin{array}{lcl}
\p_t^2 v - \p_x^2 v= 0 &\quad \mbox{on}& \R^2,\\
v(\cdot,0) = a,\ \p_tv(\cdot,0) = b &\quad \mbox{on}& \R.
\end{array}
\end{equation}
By d'Alembert's formula, it is given by
\begin{equation}\label{eq:linearwavesolution}
  v(x,t) = \frac{1}{2}\big(a(x-t) + a(x+t)\big) + \frac{1}{2}\int_{x-t}^{x+t}b(y) \dd y.
\end{equation}
Its singular support $S$ consists of the characteristic lines with slope $\pm1$ emanating from the points of
support of $a,b$, hence has Lebesgue measure zero. Denote by $\eins_{S\comp}\,v$ the almost
everywhere defined, measurable function equal to $v(x,t)$ for $(x,t) \not\in S$.
In fact, the function $\eins_{S\comp}\,v$ is piecewise constant and arises from the measure part in $b$
solely. For example, when $b(x) = \delta(x)$, we have that $\eins_{S\comp}\,v(x,t) = 1/2$ for
$|x| < t$ and $\eins_{S\comp}\,v(x,t) = 0$ otherwise.
If the function $f$ is continuous, then the nonlinear equation
\begin{equation}\label{eq:classicaldeltasemilinear}
\begin{array}{lcl}
\p_t^2 w - \p_x^2 w= f(\eins_{S\comp}\,v + w) &\quad \mbox{on}& \R^2,\\
w(\cdot,0) = 0,\ \p_tw(\cdot,0) = 0 &\quad \mbox{on}& \R
\end{array}
\end{equation}
has a unique weak solution $w \in \cC(\R^{2})$. Finally, let $u\in\cG(\R^2)$ be
the Colombeau solution to the nonlinear problem
\begin{equation}\label{eq:Colombeaudeltasemilinear}
\begin{array}{lcl}
\p_t^2 u - \p_x^2 u= f(u) &\quad \mbox{on}& \R^2,\\
u(\cdot,0) = \iota(a),\ \p_tu(\cdot,0) = \iota(b) &\quad \mbox{on}& \R.
\end{array}
\end{equation}
As we shall see shortly, $v + w$ is a delta wave, i.e. the associated distribution corresponding to the
Colombeau solution to (\ref{eq:Colombeaudeltasemilinear}).
\begin{prop}
Assume that the data $a, b \in \cD'(\R)$ have discrete support, $f$ is smooth, bounded and
globally Lipschitz. Then the generalized solution $u \in \cG(\R^{2})$ is associated
with the sum of the distributions $v \in\cD'(\R^2)$ and the continuous function $w \in \cC(\R^{2})$,
that is, $u \approx \iota(v+w)$.
\end{prop}
{\em Proof:}
The solution $u\in\cG(\R^2)$ has a representative which satisfies
\begin{equation*}
\begin{array}{lcl}
\p_t^2 u_\eps - \p_x^2 u_\eps = f(u_\eps)  &\quad \mbox{on}& \R^2,\\
u_\eps(\cdot,0) = a_\eps,\ \p_tu_\eps(\cdot,0) = b_\eps &\quad \mbox{on}& \R
\end{array}
\end{equation*}
where $a_\eps = a\ast\vphi_\eps, b_\eps = b\ast\vphi_\eps$. Further, let $v_\varepsilon$ be the classical smooth
solution to the linear wave equation (\ref{eq:classicaldeltalinear}) with data $a_\eps, b_\eps$. We have
\begin{equation}\label{eq:differencedeltasemilinear}
\begin{array}{lcl}
(\p_t^2 u - \p_x^2) (u_\eps - v_\eps - w)= f(u_\eps) - f(\eins_{S\comp}\,v + w)&\quad \mbox{on}& \R^2,\\
(u_\eps - v_\eps - w)(\cdot,0) = 0,\ \p_t(u_\eps - v_\eps - w)(\cdot,0) = 0 &\quad \mbox{on}& \R.
\end{array}
\end{equation}
We rewrite the right hand side of the first line in (\ref{eq:differencedeltasemilinear}) as
\[
f(u_\eps) - f(v_\eps + w) + f(v_\eps + w) - f(\eins_{S\comp}\,v + w).
\]
Using (\ref{eq:Linftyestimate}) with the $\rL^1$-norms, we arrive at
\begin{eqnarray*}
    \|u_\eps - v_\eps - w\|_{\rL^1(K_t)} &\leq&  T \int_0^t \|f'\|_{\rL^\infty(\R)} \|u_\eps - v_\eps - w\|_{\rL^1(K_s)}\dd s\\
          & + & T^2\|f(v_\eps + w) - f(\eins_{S\comp}\,v + w)\|_{\rL^1(K_T)}
\end{eqnarray*}
valid for $0\leq t \leq T$. Since $f(v_\eps + w) - f(\eins_{S\comp}\,v + w)$ converges to zero almost
everywhere and remains bounded, Lebesgue's theorem shows that its $\rL^1$-norm on $K_T$ converges to zero.
By Gronwall's lemma, it follows that the $L^1$-norm of
$u_\eps - v_\eps -w$ converges to zero on any $K_T$ as well. Hence $u_\eps$
converges to $v+w$ weakly, which translates into the claimed association
result. \qed

This result is an example of a wide range of much more general results on the existence of delta waves.
We refer to the papers \cite{MO86, MO91, MOWang94a, RauReed87}. In case the driving term $h$
is white noise, the generalized solution to the semilinear (stochastic)
wave equation is associated with the solution of a linear stochastic wave
equation in many cases. This has been shown e.g.\ in \cite{AlbHabRus96, MORus98a, MORus01}.
White noise or positive noise in the initial data has been considered in
\cite{MORus98b,Raj04}.
%
%
\section{Linear regularity theory}
\label{Sec:linregu}
%
%
Consider the linear wave equation (\ref{eq:classicaldeltalinear}) with initial
data distributions $a,b\in\cD'(\R)$. The classical propagation of singularities results
says that the singular support of the solution $v\in\cD'(\R^2)$ is contained in the
union of characteristic lines emanating from the singular support of the initial data.
For example, when the singular support of the initial data is situated at the origin,
the solution is a smooth function outside the forward and backward light-cone
$S = \{(x,t)\in\R^2:|x| = |t|\}$. In the Colombeau setting, this is not so, as shown by the
following example.
\begin{ex} \label{ex:deltasq}
Let $u\in\cG(\R^2)$ be the Colombeau solution to the linear wave equation
\begin{equation*}
\begin{array}{lcl}
\p_t^2 u - \p_x^2 u= 0 &\quad \mbox{on}& \R^2,\\
u(\cdot,0) = 0,\ \p_tu(\cdot,0) = \iota(\delta)^2 &\quad \mbox{on}& \R.
\end{array}
\end{equation*}
with initial data given by the square of the Dirac measure in $\cG(\R)$.
Then $u$ does not coincide with a $\Cinf$-function in the interior of the light-cone
$\{(x,t)\in\R^2:|x| < |t|\}$.
\end{ex}
In fact, by d'Alembert's formula, the solution $u \in \cG(\R^2)$
has the value
\[
   \frac{1}{2} \int_{-\infty}^\infty \iota(\delta)^2(x) \dd x
   \ = \ \mbox{class\ of\ } \left(\frac{1}{2} \int_{-\infty}^\infty \varphi_\eps^2(x) \dd x\right)_{\eps\in(0,1]}
\]
inside the forward light-cone $\{(x,t)\in\R^2:|x| < t\}$, $t>0$,
which is a generalized constant not belonging to $\C$.

The example shows that the space $\Cinf$ is not suitable for measuring regularity in Colombeau algebras.
Its role is taken over by the algebra $\Ginf$ defined in the Introduction. We have the
following result.
\begin{prop}\label{prop:linGinf}
Let $v\in\cG(\R^2)$ be the Colombeau solution to the linear wave equation
(\ref{eq:classicaldeltalinear}) with initial data $a,b,\in\cG(\R)$.
If $a,b\in \Ginf(\R\setminus\{0\})$, then $u\in \Ginf(\R\setminus S)$.
\end{prop}
{\em Proof:} We take a cut-off function $\chi\in\Cinf(\R)$ which is equal to one in a small
neighborhood of the origin and zero otherwise. Then $v = v_1 + v_2$ where $v_1\in\cG(\R^2)$
is the solution with initial data $a_1=(1-\chi)a$, $b_1=(1-\chi)b$ and $v_2\in\cG(\R^2)$
is the solution with initial data $a_2=\chi a$, $b_2=\chi b$. Using d'Alembert's formula,
we have that
\[
   \p_xv_{1\eps} = \frac{1}{2}\big(a_{1\eps}'(x+t) + a_{1\eps}'(x-t) + b_{1\eps}(x+t) - b_{1\eps}(x-t)\big)
\]
and similarly for $\p_tv_{1\eps}$. Thus if all the derivatives of $a_\eps$ and $b_\eps$ are of the same
local order, so are the derivatives of $v_{1\eps}$. This shows that $v_1$ belongs to $\Ginf(\R^2)$.
Similarly, differentiating $v_2$ once shows that the first and higher derivatives of $v_2$ vanish
outside a small neighborhood of $S$, determined by the diameter of the support of $\chi$. Since $\chi$
was arbitrary, this shows that $v_2$ is a generalized constant on $\R\setminus S$, thus belongs to
$\Ginf(\R\setminus S)$ as well. \qed

It is clear that Prop.\;\ref{prop:linGinf} can be generalized to initial data with discrete
$\Ginf$-singular support. It also holds in any space dimension \cite{MObook}.
%
%
\section{Nonlinear regularity theory}
\label{Sec:regu}
%
%
This section addresses propagation of singularities for the semilinear wave equation.
The initial data that are assumed to be regular outside the origin, but may have a singularity there. Following
the classical theory of propagation of jump discontinuities one may hope to prove that, in
the one-dimensional case, the solution will be regular (in the appropriate sense) in the interior of
the light-cone $S$. That this behavior occurs in the $\Cinf$-category for the semilinear
wave equation in one space dimension was shown by \cite{Reed}; it does not hold in higher
space dimensions or for higher order operators \cite{Rauch}, for which {\em anomalous}
singularities may occur. The question we ask is whether this transport of regularity
into the interior of the light-cone happens in the Colombeau setting as well.
We do not have a general result yet, but we will be able to prove this transport of regularity
under conditions which resemble the situation leading to {\em delta waves}.
We shall employ the subalgebra $\Ginf$ and the space $\cL_\cG^1$ defined in the Introduction
for measuring the regularity.

We consider the semilinear wave equation (\ref{eq:semilinearwave}). As before, we assume
that the function $u \to f(u)$ is smooth with all derivatives of at most
polynomial growth as $|u|\to \infty$, and that it satisfies a global Lipschitz estimate,
i.e., has a globally bounded derivative. We take initial data of the form
\[
   u_i = r_i + s_i \in \cG(\R),\quad i = 0,1
\]
where
\[
   r_i \in \cL_\cG^1(\R) \ {\rm and}\ \supp (s_i) = \{0\},\quad i = 0,1.
\]
Define the generalized complex number $M\in\widetilde{\C}$ as the class of
\[
   M_\eps = \frac{1}{2} \int_{-\infty}^\infty s_{1\eps}(x) \dd x
\]
where $(s_{1\eps})_{\eps\in (0,1]}$ is a representative of $s_1$. We shall write
$M \approx m$ if $M_\eps$ converges to a complex number $m\in\C$ and $|M| \approx \infty$
if $|M_\eps| \to \infty$ as $\eps\to 0$. Again $S = \{(x,t)\in\R^2:|x| = |t|\}$ denotes
the forward and backward light-cone.
\begin{prop} \label{prop:deltareg}
Assume that the function $f$ is smooth, globally Lipschitz, and all its derivatives are
polynomially bounded.  Let $u\in \cG(\R^2)$ be the solution to the semilinear wave equation
(\ref{eq:semilinearwave}) with initial data $u_0, u_1$ as described above. If either

(a) $f$ is globally bounded and $M \approx m$ for some $m\in\C$, or

(b) $f$ is globally bounded, $\lim_{|y|\to\infty} f(y)$ exists and $|M| \approx \infty$

then $u\in \Ginf\big(\R^2\setminus S\big) +  \cL_\cG^1\big(\R^2)$.
\end{prop}
{\em Proof:}
We first derive the estimates for $(x,t)$ in the upper half-plane $\hp$.
Denote by $\eins_\Sigma$ the characteristic function of the solid light-cone
$\Sigma = \{(x,t):t\geq 0,|x|\leq t\}$. Let $(r_{i\eps})_{\eps\in (0,1]}$, $(s_{i\eps})_{\eps\in (0,1]}$,
$i = 0,1,$ be representatives of the initial data, where we may assume that $\supp(s_{i\eps}) \subset [-\eta,\eta]$
with $\eta$ as small as we wish. Let $u_\eps,v_\eps \in\Cinf(\hp)$, $w_\eps \in \rL_{\mathrm{loc}}^1(\hp)$ be the solutions to
\begin{equation*}
\begin{array}{lll}
\left(\p_t^2 - \p_x^2\right)u_\eps =  f(u_\eps), &\quad u_\eps(\cdot,0) = u_{0\eps},&
                      \p_tu_\eps(\cdot,0) = u_{1\eps},\vspace{2pt}\\
\left(\p_t^2 - \p_x^2\right)v_\eps =  0,  &\quad v_\eps(\cdot,0) = s_{0\eps}, &
                      \p_tv_\eps(\cdot,0) = s_{1\eps} \vspace{2pt},\\
\left(\p_t^2 - \p_x^2\right)w_\eps =  f(M_\eps\eins_\Sigma + w_\eps), &\quad w_\eps(\cdot,0) = r_{0\eps},&
                      \p_tw_\eps(\cdot,0) = r_{1\eps}.
\end{array}
\end{equation*}
Then
\[
  \left(\p_t^2 - \p_x^2\right)\big(u_\eps - v_\eps - w_\eps\big)
    = f(u_\eps) - f(v_\eps + w_\eps) + f(v_\eps + w_\eps) - f(M_\eps\eins_\Sigma + w_\eps)
\]
with zero initial data. The difference of the first two terms on the right hand side can be estimated by
the $\rL^\infty$-norm of $f'$ times $|u_\eps - v_\eps - w_\eps|$, while the difference of the last two terms
vanishes off an $\eta$-neighborhood of the light-cone $S$, because $v_\eps(x,t) = M_\eps$ for
$|x| < t-\eta$. Taking a trapezoidal region $K_T$ as in (\ref{eq:ItKT}), the boundedness of $f$,
inequality (\ref{eq:Linftyestimate}) and Gronwall's lemma give an estimate of the form
\begin{equation} \label{eq:conv1}
   \|u_\eps - v_\eps - w_\eps\|_{\rL^1(K_T)} \leq C\eta
\end{equation}
for some constant $C>0$ and all $\eps \in (0,1]$. In the case (a), let $w \in \rL_{\mathrm{loc}}^1(\hp)$
be the solution to
\begin{equation*}
\left(\p_t^2 - \p_x^2\right)w =  f(m\eins_\Sigma + w), \qquad w(\cdot,0) = r_0,\quad \p_tw(\cdot,0) = r_1
\end{equation*}
where $r_i$ is the limit in $\rL_{\mathrm{loc}}^1(\R)$ of $r_{i\eps}$ as $\eps\to 0$, $i = 0,1$.
Now
\[
  \left(\p_t^2 - \p_x^2\right)\big(w_\eps - w)
  = f(M_\eps\eins_\Sigma + w_\eps) - f(M_\eps\eins_\Sigma + w) + f(M_\eps\eins_\Sigma + w) - f(m\eins_\Sigma + w)
\]
with initial data given by $r_{i\eps} - r_i$, $i = 0,1$. The difference of the first two terms on the
right hand side is bounded by the $\rL^\infty$-norm of $f'$ times $|w_\eps - w|$. By assumption,
the difference of the last two terms converges to zero almost everywhere. Using inequality
(\ref{eq:Linftyestimate}) with $\rL^1$-norms in place of the $\rL^\infty$-norms and Gronwall's lemma
as above shows that
\begin{equation} \label{eq:conv2}
  \|w_\eps - w\|_{\rL^1(K_T)} \to 0 \quad {\rm as}\ \eps\to 0.
\end{equation}
In case (b), defining $w \in \rL_{\mathrm{loc}}^1(\hp)$ as the solution to
\begin{equation*}
\left(\p_t^2 - \p_x^2\right)w =  L\eins_\Sigma, \qquad w(\cdot,0) = r_0,\quad \p_tw(\cdot,0) = r_1,
\end{equation*}
the same argument as above leads to the convergence result (\ref{eq:conv2}) in this case as well.
Combining (\ref{eq:conv1}) which holds for arbitrarily chosen $\eta >0$ with (\ref{eq:conv2}) shows that
$u_\eps - v_\eps$ converges to $w$ in $\rL^1(K_T)$ as $\eps\to 0$.
By Prop.\;\ref{prop:linGinf}, $v_\eps$ enjoys the $\Ginf$-estimate (\ref{eq:regufu}) off the light-cone $S$.
Thus $u_\eps = v_\eps + (u_\eps - v_\eps)$ defines an element of
$\Ginf\big((\hp)\setminus S\big) +  \cL_\cG^1\big(\hp\big)$.

The arguments for $(x,t)$ in the lower half-plane $\R\times (-\infty, 0]$ are the same with $-M$ in place
of $M$.
\qed

\begin{rem}
The hypotheses on the generalized constant $M$ in Prop.\;\ref{prop:deltareg} are satisfied when the
term $s_1$ in the initial data is a polynomial in the Dirac measure and its derivatives. In fact,
when $s_1 = \p^\alpha\iota(\delta)$, we have $M\approx 1$ for $\alpha = 0$ and $M\approx 0$ for $\alpha > 0$.
When $s_1 = \pi(\iota(\delta))$ for some polynomial function $\pi$, only the cases $M\approx m$ for some
$m\in\C$ or $|M|\approx\infty$ can occur.
\end{rem}

Prop.\;\ref{prop:deltareg} shows that regularity of the type $\Ginf + \cL_\cG^1$ is propagated into the region inside
the light-cone. It is clear that it can be generalized in various ways: for example, the support of the singular part
$s_0, s_1$ of the data could consist of a discrete set rather than just a point. However, it should be noted
that continuous dependence of the regularized solutions on the data in terms the $\rL^1$-norm enters into the
argument, and such a property depends decisively on the particular equation, the space dimension and the
nonlinearity $f$. Prop.\;\ref{prop:deltareg} exploits such special properties and thus falls short of providing a
prototypical description of nonlinear propagation of regularity for Colombeau solutions, which remains a
challenging open question.

Another promising approach to measuring regularity in the nonlinear case is
the Colombeau-H\"older-Zygmund-scale which has been introduced and applied to
nonlinear scalar first order equations in \cite{Hoer04}.

\end{document}